\documentclass{article}
\usepackage[utf8]{inputenc}
\usepackage[left=1cm,right=1cm,top=2cm,bottom=2cm]{geometry}

\pdfoutput=1
\usepackage{lineno,hyperref}
\modulolinenumbers[5]
\usepackage{amsmath}

\usepackage{amssymb}
\usepackage{float}
\usepackage{mathtools}
\usepackage{lipsum}
\usepackage{caption}
\usepackage{subcaption}
\usepackage{commath}
\usepackage{cases}

\usepackage{amsthm}
\usepackage{graphicx}
\usepackage{listings}
\usepackage{multicol}

\newcommand{\pt}{{\partial}}
\newcommand{\nn}{{\mathbf{n}}}

\newcommand{\etaeta}{{\bf \eta}}
\newcommand{\uu}{{\bf u}}

\begin{document}
 \title{\noindent A mathematical model for meniscus cartilage regeneration}
\maketitle
\normalsize
\begin{center}
\author{Elise Grosjean \footnotemark[1],}
\author{Bernd Simeon \footnotemark[1],}
\author{Christina Surulescu \footnotemark[1]}
\end{center}

\footnotetext[1]{Felix-Klein-Institut für Mathematik, RPTU Kaiserslautern-Landau, Kaiserslautern, 67663, Deutschland}

\date

\abstract{We propose a continuous model for meniscus cartilage regeneration triggered by two populations of cells migrating and (de)differentiating within an artificial scaffold with a known structure. The described biological processes are influenced by a fluid flow and therewith induced deformations of the scaffold. Numerical simulations are done for the corresponding  dynamics within a bioreactor which was designed for performing the biological experiments.}

\maketitle

\section{Introduction}\label{sec1}
A meniscus can be defined as a crescent-shaped fibro-cartilaginous tissue that is responsible for the structural integrity of the knee.  Meniscus tears are among the most common injuries in persons doing contact sports. They can lead to severe cartilage degeneration and therewith associated with severe chronic knee pain, stiffness and immobility. Numerous patients finally have to face a total knee replacement. There is a recent trend in  orthopedic surgery to promote healing by regeneration and repair rather than replacing the damaged tissue, in order to avoid late degeneration. Accordingly, there is a strong need for a better understanding of the biological processes inside the
meniscus tissue, and regenerative meniscus substitutes are increasingly coming into focus.
In such substitutes, living cells are combined with suitable biochemical and physicochemical factors and subjected to engineering methods. Often, such scaffolds possess a nonwoven-type structure that mimics the biological tissue.
For proper functioning, the generated tissues require certain mechanical and structural properties.\\

Whereas meniscus tissue engineering has recently attracted much interest, mathematical models for the dynamics of  involved phenomena are relatively scarce. Hallmark issues of meniscus regeneration relate to degradation of engineered fibers, migration/differentiation of stem cells into/within the scaffold, and production of tissue by  chondrocytes. Several factors are thereby believed to play an essential role, prominently the stem cell (de)differentiation triggered by mechanical stress \cite{B6:1}, tissue stiffness \cite{Mosavi}, topography of the scaffold \cite{Ghasemi}, or by chemical cues present in the extracellular space \cite{B6:11}.
We are not aware of any continuous settings addressing meniscus repair from the said perspective; \cite{geris-etal} proposed a pure macroscopic model in a related context, accounting for biochemical, but not for mechanical influences and not being able to capture the topography of underlying tissue.
Most continuum approaches feature reaction-diffusion-(transport) equations (RD(T)Es) for
the evolution of macroscopic cell 
densities interacting with chemical cues (chemotaxis) and/or tissue (haptotaxis). 
Yet other macroscopic formulations use a multiphase approach where the cell populations are seen as components of a mixture also containing fluid(s) and/or tissue, possibly also soluble matter (acting as chemical cues), and which relies on mass and momentum balance for each of the involved phases, supplemented with appropriate closure laws, see e.g. \cite{B6:2,LemonKing07}; the review \cite{Klika} of multiphase cartilage mechanical modeling explicitly excludes descriptions of cell behavior involved in the process. 
A connection between multiphase models and RDTE systems for (tumor) cell migration and spread in the extracellular matrix has been proved for any space dimension in \cite{KSZ22}.

This work is organized as follows: in  Section \ref{sec:model} we state the PDE-ODE system for the dynamics of two cell populations involved in cartilage expression, along with those of the soluble and insoluble  signals in the extracellular space.  Section \ref{sec:numerics} is concerned with numerical simulations of this system and its extensions to a model also accounting for the fluid effects and the deformations of the meniscus, which is seen as a poroelastic material. A coupling strategy for the mechanical and the biological building blocks of our model is furthermore presented, along with numerical results in Subsection \ref{sec:coupling}.

\section{Model setup}\label{sec:model}

We consider a first model for the dynamics of adipose derived stem cells (ADSCs) differentiating into and interacting with chondrocytes, under chemical and mechanical environmental influences. Chondrocytes produce cartilage and hyaluron, also uptaking the latter. The cells migrate, differentiate, and proliferate inside an artificial scaffold with given topology, which does not infer resorbtion by the cells and whose fibers are impregnated with hyaluron, that acts as a (nondiffusing) chemoattractant for the ADSCs: 
\begin{align}
\partial_t c_1 &= a_1 \Delta c_1 - \nabla \cdot (b_1 c_1 \nabla h) - \nabla \cdot (b_2 c_1 \nabla k) - \alpha_1(k,S) c_1 +  \alpha_2(k, S) c_2 +\beta c_1\left(1 -c_1-c_2-k \right), 	\label{omodel-c_1} \\
\partial_t c_2  &=  \Delta c_2 + \alpha_1(k, S) c_1 - \alpha_2(k,S)c_2,\label{omodel_c_2}                \\
\partial_t h &= - \gamma_1 h c_2 + \dfrac{c_2}{1+c_2},\label{omodel-h}      \\
\partial_t k&=  - \delta_1 c_1 k + c_2,   	\label{omodel-tau}                                         
\end{align}
where $c_1$ is the density of ADSCs which differentiate into and interact with chondrocytes $c_2$ producing cartilage $k$ and $h$ is the density of hyaluron. All parameters are assumed
to be positive; $\alpha_1(k, S)$ and $\alpha_2(k, S)$ are uniformly bounded functions of $k$ and $S$, where $S$ represents a quantity related to the mechanical stress exerted on the cells. The above equations are already nondimensionalized. The domain $\Omega_p \subset \mathbb{R}^d,\ d \leq 3$ is bounded, with a smooth boundary. The boundary conditions are of the zero-flux type
\begin{equation}
	\label{bco}
	-\dfrac{\partial c_1}{\partial \nu} + b_1 c_1 \dfrac{\partial h}{\partial \nu} + b_2 c_1 \dfrac{\partial k}{\partial \nu} = \dfrac{\partial c_2}{\partial \nu} = 0 ~~\text{on}~~ \partial \Omega_p \times (0, T),
\end{equation}
where $\nu$ is the outward unit normal on the boundary of $\Omega_p$. The initial conditions are:
\begin{equation}
	\label{oic}
	c_1(x, 0) = c_{10}(x) , ~~ c_2(x, 0) = c_{20}(x), ~~ h(x, 0) = h_0(x), ~~ k(x, 0) = k_0(x), ~~ x \in \Omega_p.
\end{equation}
\noindent
This macroscopic formulation of space-time cell population dynamics provides a preliminary, phenomenological description of putatively essential physiological processes
related to cell proliferation, differentiation, and motility, as well as dynamics of tissue. It is in fact a simplified version of the macroscopic system 
\begin{align}
&\partial_tc_{1}-\nabla\nabla :\left (\mathbb D_1c_{1}\right )+\nabla \cdot \left (\chi (B(h,k))\mathbb D_1\nabla B(h,k)c_{1}\right )=-\alpha_1(k,S)c_{1}
+\alpha_2(k,S)\frac{\omega_1}{\omega_2}c_{2}+\beta c_{1}\left(1-c_1-c_2\right )\label{macro-c1}\\
\vspace{-0.25cm}
&\partial_tc_{2}-\nabla\nabla :\left (\mathbb D_2c_{2}\right )=\alpha_1(k,S)\frac{\omega_2}{\omega_1}c_{1}-\alpha_2(k,S)c_{2}\label{macro-c2}
\end{align}
supplemented with the two macroscopic ODEs 	\eqref{omodel-h}, \eqref{omodel-tau} for $h$ and $k $. 
Thereby, $\nabla\nabla :(\mathbb D\mathbb c)=\nabla \cdot (\nabla \cdot \mathbb Dc+\mathbb D\nabla c)$ represents myopic diffusion, 
whereby the diffusion coefficient $\mathbb D$ encodes the orientation distribution of scaffold fibers.  The quantity $B(h,k)$ is an affine function of $h$ and $k $.  
\noindent
System \eqref{macro-c1}, \eqref{macro-c2} can be obtained upon starting from lower scales 
and performing a parabolic upscaling, similarly to e.g., \cite{CDKS,CS-21}. More details will be available in a forthcoming work.
\section{Numerical methods} \label{sec:numerics}
At the macroscopic scale, the meniscus tissue
is a poroelastic medium that can be modeled by Biot's equations
\begin{subequations}
	\begin{align}
\rho_s \pt_{tt} \boldsymbol{\eta}_p - \nabla\cdot {\boldsymbol{\sigma}_p}(\boldsymbol{\eta}_p,p_p) & =  0, 
\label{eq:B61} \\
\partial_t \left( \frac{1}{M} p_p + \nabla \cdot (\alpha \boldsymbol{\eta}_p) \right)
+ \nabla \cdot \mathbf{u}_p & =  0,  \label{eq:B62}
\end{align}
\end{subequations}
for the displacement field $\boldsymbol{\eta}_p(x,t)$ and the pressure $p_p(x,t)$ in a domain
$\Omega_p \subset \mathbb{R}^d$, with additional boundary and initial conditions.
Here, $\rho_s$ stands for the solid phase density while $M$ and $\alpha$ are Biot's modulus and coefficient, respectively. The stress tensor $\boldsymbol{\sigma}_p$ is given by an appropriate constitutive equation and the fluid flux satisfies Darcy's law 
\vspace{-0.3cm}
\begin{equation}
\label{eq:darcyeq}
\mathbf{u}_p = -\mathbf{K}(\nabla p - \rho_f \mathbf{g})/\mu,
\end{equation} with permeability matrix $\mathbf{K}$, fluid phase density $\rho_f$, and  viscosity $\mu$. Given suitable material properties and geometry data, 
model  \eqref{eq:B61}-\eqref{eq:B62}  can be solved by the Finite Element Method (FEM),
as demonstrated in \cite{B6:6}. 
For the active biological processes inside the tissue, however,
a more detailed description is required, that takes the temporal
change of the solid and fluid phases as well as the material properties of collagen gel into account. In \cite{B6:2}, such a model is proposed based on fractional volumes for solid and fluid
phases and a Maxwell-type constitutive equation.
Aside from the basic equations of mass and momentum, an additional evolution equation expresses the spreading of cells into the structure, as given by system \eqref{omodel-c_1}-\eqref{omodel-tau}. Overall, this results in a complex coupled problem in which the cell densities appear as additional unknowns. The goal of this section is to propose a loosely coupling strategy to connect the Biot-Darcy system to the cell densities evolution and we therefore present 3D results ($d=3$) on both models, which have been implemented using FreeFem++ \cite{hecht2012new} in parallel with PETSc \cite{petsc}.
\subsection{Simulations of the macroscopic advection-diffusion-reaction equations}
\label{macroeq}
We proceed with direct simulations of the model \eqref{omodel-c_1}-\eqref{omodel-tau} on the scaffold, denoted $\Omega_p$. Its numerical scheme should be locally mass conservative, thus we have decided to employ a first order Non-symmetric Interior Penalty discontinuous Galerkin (NIP dG) scheme in space  \cite{di2011mathematical}. 
We define a mesh $\mathcal{T}_h$ of $\Omega_p$ and seek solutions $c_1$ and $c_2$ in the broken polynomial space $\mathbb{P}^1_d( \mathcal{T}_h)$ given by $\mathbb{P}_d^1(\mathcal{T}_h):=\{ u \in L^2(\Omega_p)\ | \ \forall T \in \mathcal{T}_h, v_{|_T} \in \mathbb{P}_d^1(T) \}$, ~whereas we are looking for $h$ and $k$ on the classical $\mathbb{P}^1_d( \Omega_p)$ FE space.  We have used 88 634 degrees of freedom, and a time step $\Delta t=0.1$.
Multiplying by test functions $(\nu_{c1}, \nu_{c2}, \nu_{h},\nu_{k})$ and integrating over $\Omega_p$, system \eqref{omodel-c_1}-\eqref{omodel-tau} becomes
\begin{equation}
  \begin{cases}
  \label{NIPsys1}
    & (\pt_t c_1, \nu_{c1}) + ( a_1 \nabla c_{1},\nabla \nu_{c1}) + ([c_1],\{a_1 \nabla \nu_{c1}\})_{\Gamma} - ([\nu_{c1}],\{a_1 \nabla c_1\})_{\Gamma} \\
  & \quad \quad \quad  - (\textbf{v}\ c_1,\nabla \nu_{c1}) + ((\textbf{v} c_1)^{\uparrow},[\nu_{c1}])_{\pt \Omega_p}   + (\alpha_1 c_1 -\alpha_2 c_2  - \beta c_1(1-c_1-c_2-k), \nu_{c1}) + (\eta[c_{1}],[\nu_{c1}])_{\Gamma}=0,  \\
&(\pt_t c_2, \nu_{c2}) + (  \nabla c_{2},\nabla \nu_{c2}) + ([c_2],\{\nabla \nu_{c2} \})_{\Gamma}- ([\nu_{c2}],\{ \nabla c_2\})_{\Gamma} - (\alpha_1 c_1 -\alpha_2 c_2, \nu_{c2}) +(\eta[c_{2}],[\nu_{c2}])_{\Gamma}= 0,  \\
& (\pt_t h, \nu_h) + (\gamma_1\ h \ c_2,\nu_h) - ( \frac{c_2}{1+c_2},\nu_h)=0, \\
& (\pt_t k, \nu_k) + (\delta_1\ k \ c_1,\nu_k) - ( c_2,\nu_k)=0, \\
  & c_1(0) = c_1^0, \ c_2(0) = c_2^0, \ h(0) = h_0,
\end{cases}
\end{equation}
where $\nabla$ refers to the broken gradient, $\Gamma $ represents all the interfaces of the mesh, $\eta$ is the penalization parameter, $\mathbf{v}=b_1 \nabla h+b_2 \nabla k$, $(\cdot, \cdot)$ refers to the $L^2(\Omega_p)$ inner product, $(\cdot)^{\uparrow}$ is the upwind flux, and $[\cdot]$ and $\{ \cdot \}$ refer to jumps and means. 
 The nonlinear system \eqref{NIPsys1} has then been implicitly discretized in time and solved with a Newton algorithm.
\subsection{Bioreactor simulations}
The application of mechanical loads was investigated as an important stimulus for cell growth \cite{B6:1}.
A major challenge lies in the numerically efficient coupling of the processes at the cell level with the macroscopic behavior and the mechanical properties of the tissue.
The scaffold is integrated in a 3D printed perfusion chamber which is embedded in a bioreactor. The latter enables mechanical stimulation via an alternating fluid passing through the perfusion chamber tubes, releasing the pressure. From a numerical standpoint, it can be modeled by Biot-Darcy equations \eqref{eq:B61}-\eqref{eq:B62}  in $\Omega_p$, coupled to a nonstationary Stokes problem in $\Omega_f$ (corresponding to the tubes of the perfusion chamber), where the whole spatial domain is $\Omega= \Omega_f \cup \Omega_p$. See in Figure \ref{firstplots} the domain decomposition. We have used a Nitsche's penalization approach that allows to impose the interface conditions between the free fluid part in the tubes and the tissue \cite{biotdarcystokes1} and have adapted it to our boundary conditions.
\begin{itemize}
    \item In the free fluid region denoted $\Omega_f$, we denote by $\mathbf{n}_f$ the outward unit normal vector to the boundaries $\Gamma_f=\pt \Omega \cap \Omega_f=\Gamma_I \cup \Gamma_{f,W} \cup \Gamma_{in} \cup \Gamma_{out},$ where $\Gamma_I$ represents the interface between $\Omega_f$ and $\Omega_p$. \\
    The velocity-pressure couple $(\uu_f , p_f )$ satisfies the unsteady Stokes equations 
\begin{equation}
\label{eq:unstokes}
  \rho_f \pt_t \uu_f  -  \nabla \cdot \boldsymbol{\sigma}_f(\uu_f,p_f)=0,   \quad \mathrm{ and }  \quad  \nabla \cdot \uu_f=0, \quad \mathrm{ in }\  \Omega_f,
\end{equation}
for the fluid stress tensor given by $\boldsymbol{\sigma}_f(\uu_f,p_f):=-p_f I + 2 \mu D(\uu_f),  \textrm
    {with } D(\uu_f)=\frac{1}{2} (\nabla \uu_f + \nabla \uu_f^T)$ and where  $\rho_f$ stands for the fluid phase density.  We complete \eqref{eq:unstokes} by the following boundary conditions: \\
    $\uu_f = 0, \textrm{ on } \Gamma_{f,W}$, $\boldsymbol{\sigma}_f \ \nn_f = -p_{in}(t) \ \nn_f \textrm{ and } \uu_f \times \nn_f = 0,\textrm{ on } \Gamma_{in},$ and $\boldsymbol{\sigma}_f \ \nn_f =0, \textrm{ on } \Gamma_{out}$.\vspace{-0.075cm}
\item In the tissue $\Omega_p$, $\Gamma_p= \pt \Omega \cap \Omega_p=\Gamma_{p,W} \cup \Gamma_I$. With $\lambda_p$ and $\mu_p$ the Lamé parameters, the poroelastic stress $\boldsymbol{\sigma}_p$ is defined by:
\begin{equation}
\label{eq:sigmap}\boldsymbol{\sigma}_p(\etaeta_p,p_p)=\boldsymbol{\sigma}_e(\etaeta_p) - \alpha p_p I,  \quad  \mathrm{and }  \quad \boldsymbol{\sigma}_e(\etaeta_p)=\lambda_p (\nabla \cdot \etaeta_p)\ I + 2 \mu_p D(\etaeta_p),
\end{equation}
and it is subject to \eqref{eq:B61},\eqref{eq:B62} and \eqref{eq:darcyeq} with 
 $p_p=0, \textrm{ on }   \Omega_{p,W},\ \mathrm{ and } \  \boldsymbol{\eta}_p=0,  \textrm{ on }   \Omega_{p,W}.$ 
\end{itemize}
We set first all variables to $0$ and consider a referential fixed domain. Combining with Nitsche's penalization \eqref{eq:B61}-\eqref{eq:B62}-\eqref{eq:darcyeq} and \eqref{eq:unstokes}, we arrive to a complex formulation as given in \cite{biotdarcystokes1} that we have implemented in parallel  with ($\mathcal{P}$1b, $\mathcal{P}^1$,$\mathcal{RT}0$,$\mathcal{P}^0$,$\mathcal{P}^1$)  elements for the solutions $(\uu_f,p_f,\uu_p,p_p,\etaeta_p)$. In figures \ref{firstplots} and \ref{secondplots} are illustrated results with a time step $ \Delta t= 0.1$ and the parameters presented in table \ref{tab:parametervaluesBiot}, which are obtained from literature \cite{MechaParamsLit, MechaParamsLit2}. The pressure boundary condition at the entrance of the perfusion chamber is set to 
$ p_{in}(t)=~p_{max} \ \sin(\pi \ t)$. 
\begin{table}[!ht]
    \centering
    \begin{tabular}{|c|c|c||c|c|c|}
        {Maximum fluid pressure (MPa)} & $ p_{max} $&   $10$ &    {Dynamic viscosity} (MPa.s) &$ \mu_f$& $1\times 10^{-9}$ \\
        {Poroelastic wall density} (kg/$m^3$) &$\rho_p$& $1.1 \times 10^3$  & {Permeability} ($m^4/Ns$) &$\kappa$ & 1E-14\\
        { Fluid density } (kg/$m^3$)  &$\rho_f$ &$10^{3}$ &   Initial porosity (\%) & $\Phi$ &0.8 \\
        {Young's modulus} (MPa) & E & $80 $ & Mass storativity (MPa$^{-1}$ ) & $\frac{1}{M}$ &$6.89 \times 10^{1}$  \\
        {Poisson's ratio } &$\nu$ & 0.167 &Biot–Willis constant&$ \alpha$ & $1.0$ \\
    \end{tabular}
  
    \caption{Poroelasticity and fluid parameters}
\label{tab:parametervaluesBiot}
\end{table}

\vspace{-0.5cm}
\subsection{Coupling strategy}\label{sec:coupling}
We aim at incorporating the mechanical stimulus into \eqref{NIPsys1} in order to study its impact on the densities of the cells.
 We therefore propose a loosely coupling inspired by previous studies \cite{Andreykiv_2007} where the stimulus $S$ is defined as the sum of the strain and the fluid velocity in which scaling constants are used for each stimulus. The authors showed that $S$ would encourage chondrocyte differentiation and cartilage synthesis if it fell between two values, $S_{min}$ and $S_{max}$. We thus propose a mapping between $\boldsymbol{\sigma}_p$ and the rates $\alpha_1$ and $\alpha_2$, which become dependent on time, space and on stress.

\begin{figure}[!h]
     \centering
     \begin{subfigure}[b]{0.2\textwidth}
         \centering
         \includegraphics[width=1.\textwidth]{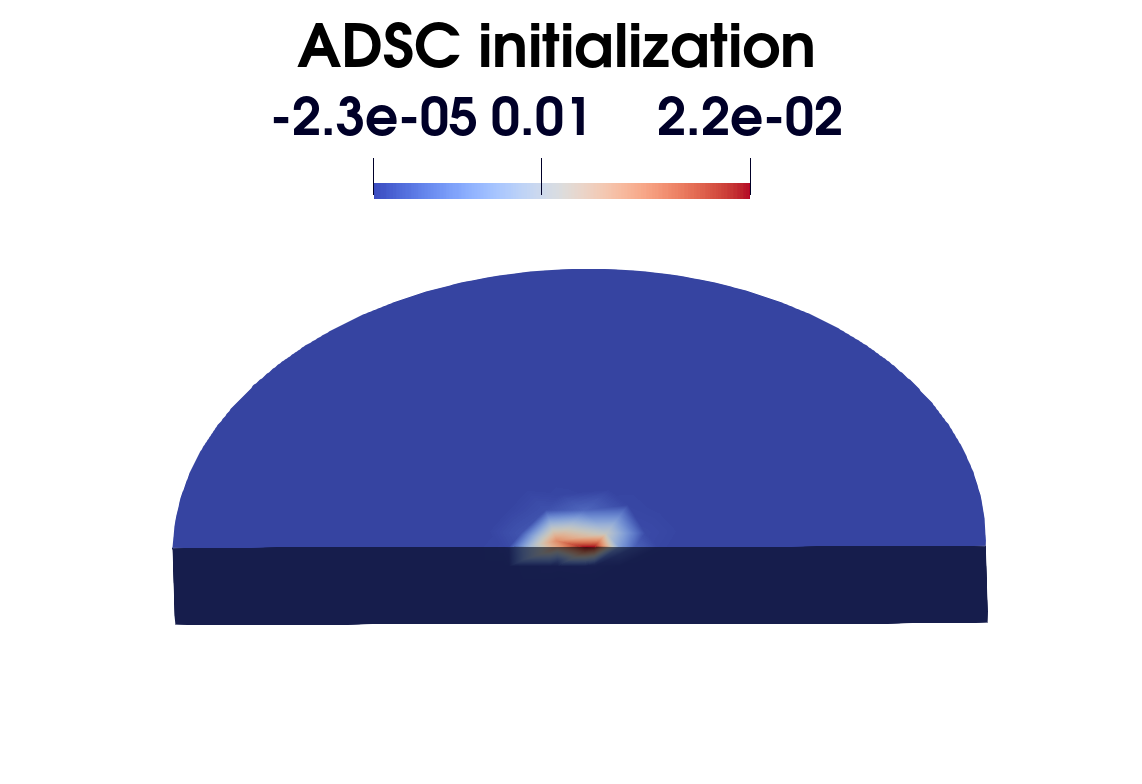}
     \end{subfigure}
          \begin{subfigure}[b]{0.2\textwidth}
         \centering
         \includegraphics[width=1.\textwidth]{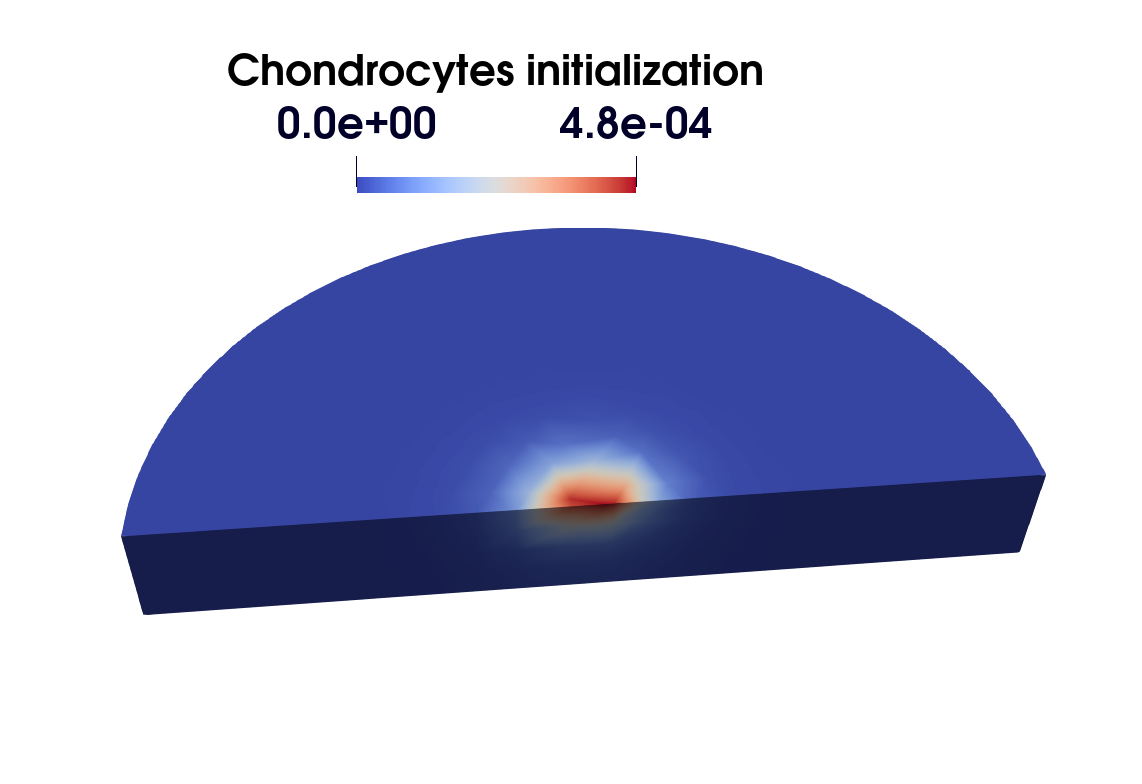}
     \end{subfigure} 
          \begin{subfigure}[b]{0.2\textwidth}
         \centering
         \includegraphics[width=1.\textwidth]{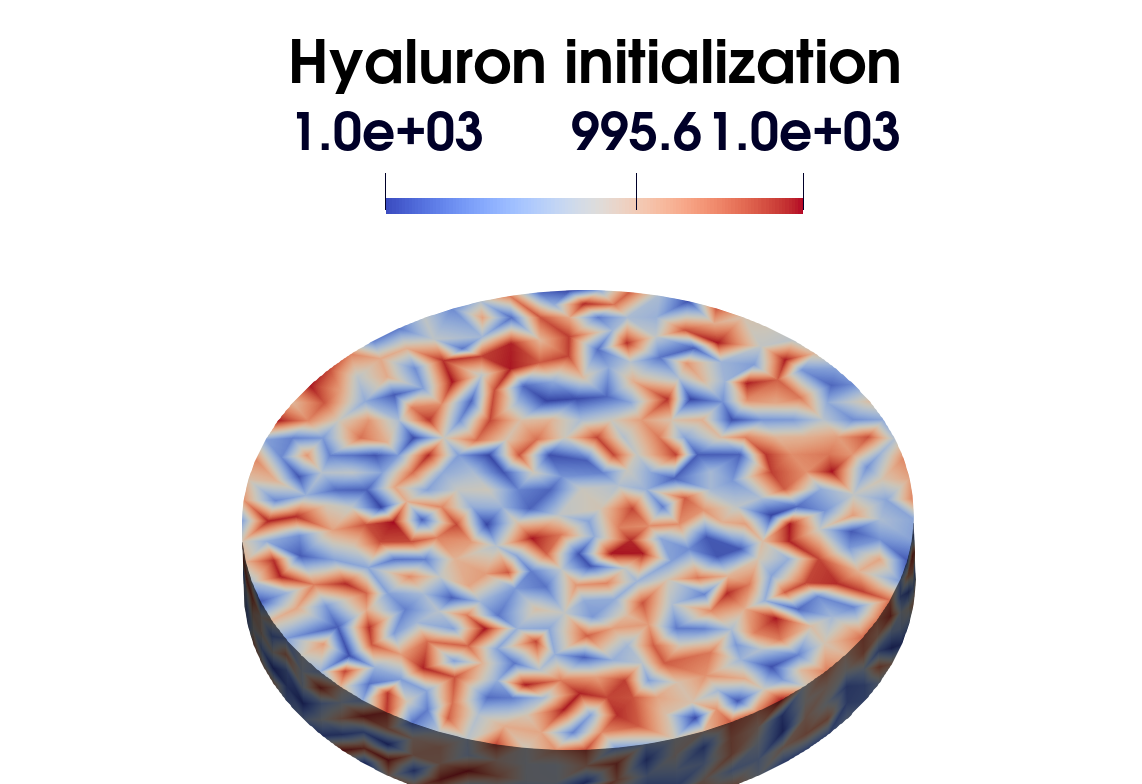}
     \end{subfigure}
          \begin{subfigure}[b]{0.2\textwidth}
         \centering
         \includegraphics[width=1.\textwidth]{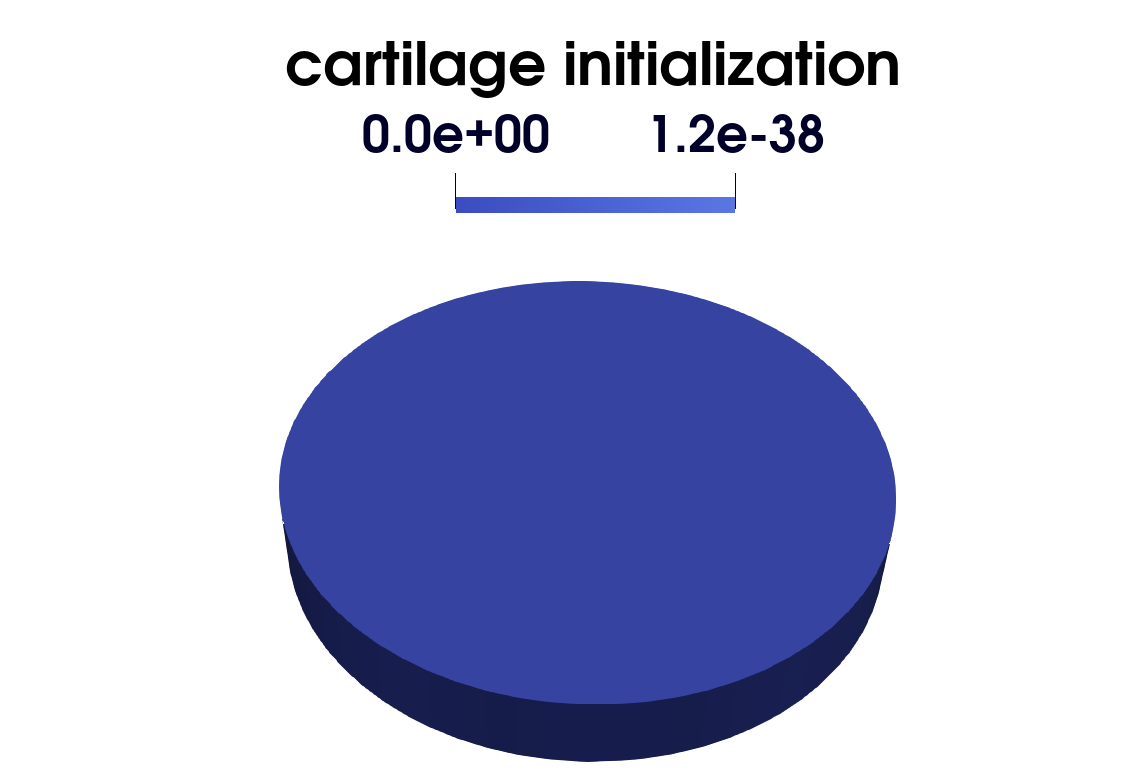}
     \end{subfigure}
            \begin{subfigure}[b]{0.18\textwidth}
         \centering
         \includegraphics[width=1.\textwidth]{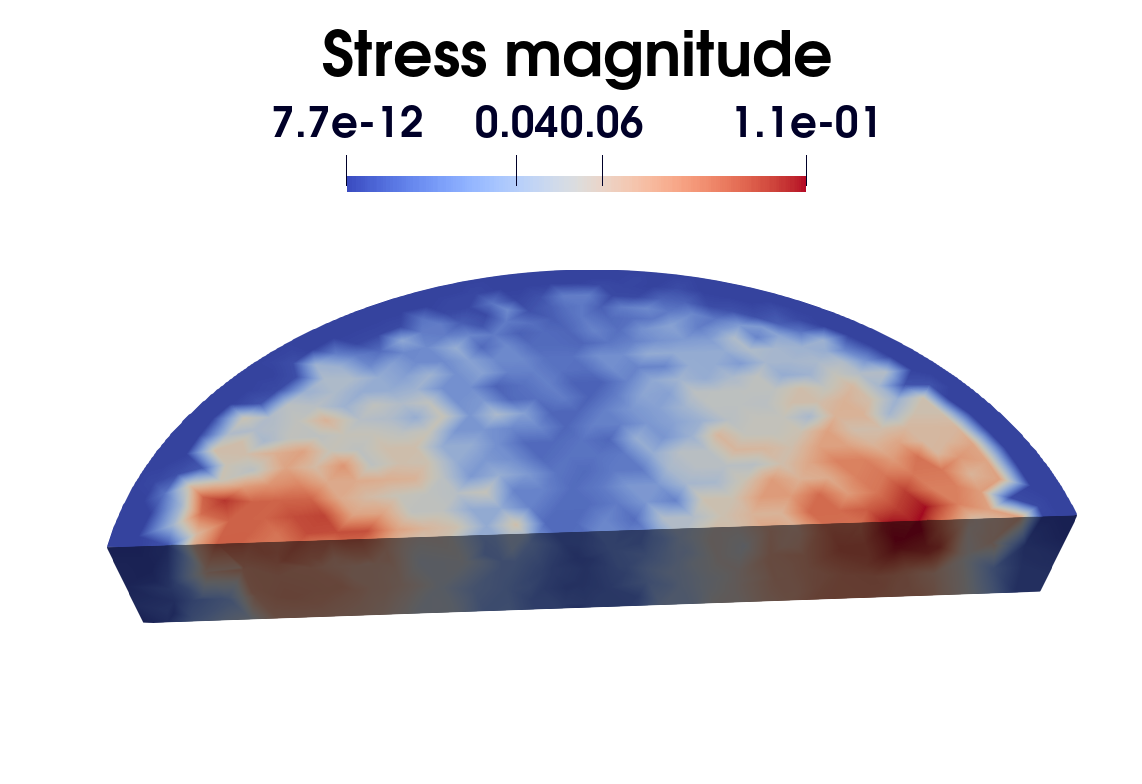}
         \end{subfigure}
     \caption{Initialization}
     \label{Cellresults}
     \end{figure}

\begin{figure}[!h]
     \centering
     \begin{subfigure}[b]{0.2\textwidth}
         \centering
         \includegraphics[width=1.\textwidth]{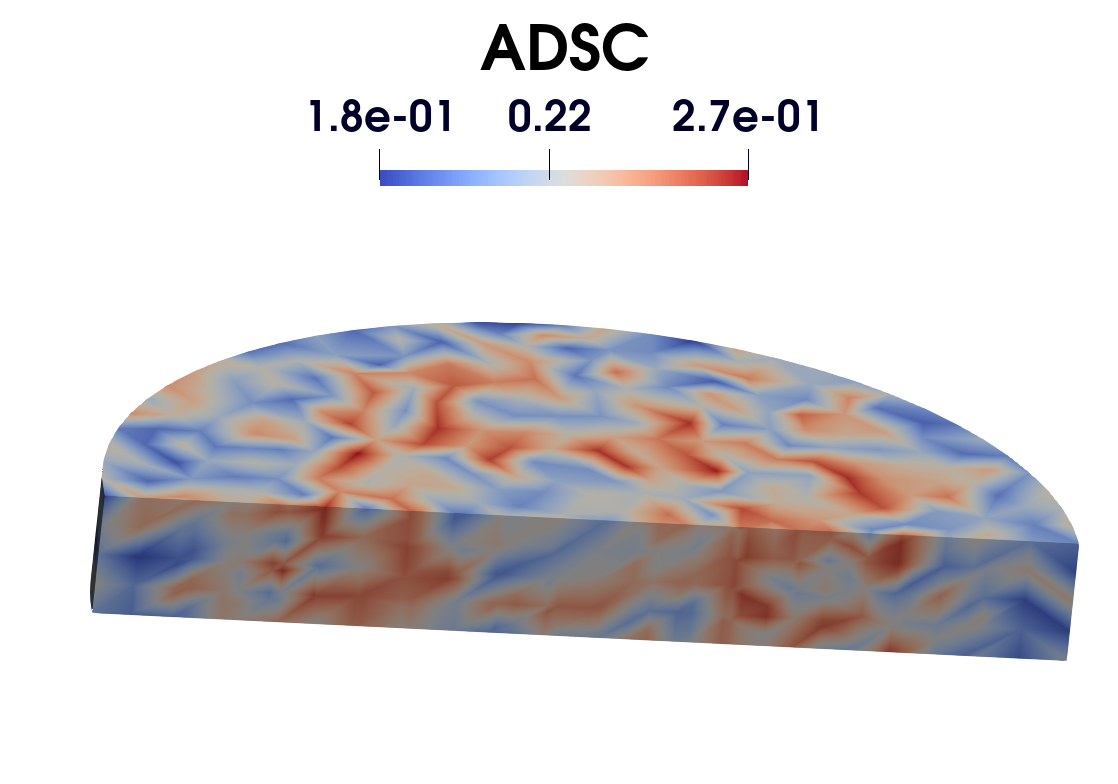} 
     \end{subfigure} \quad \quad
          \begin{subfigure}[b]{0.2\textwidth}
         \centering
         \includegraphics[width=1.\textwidth]{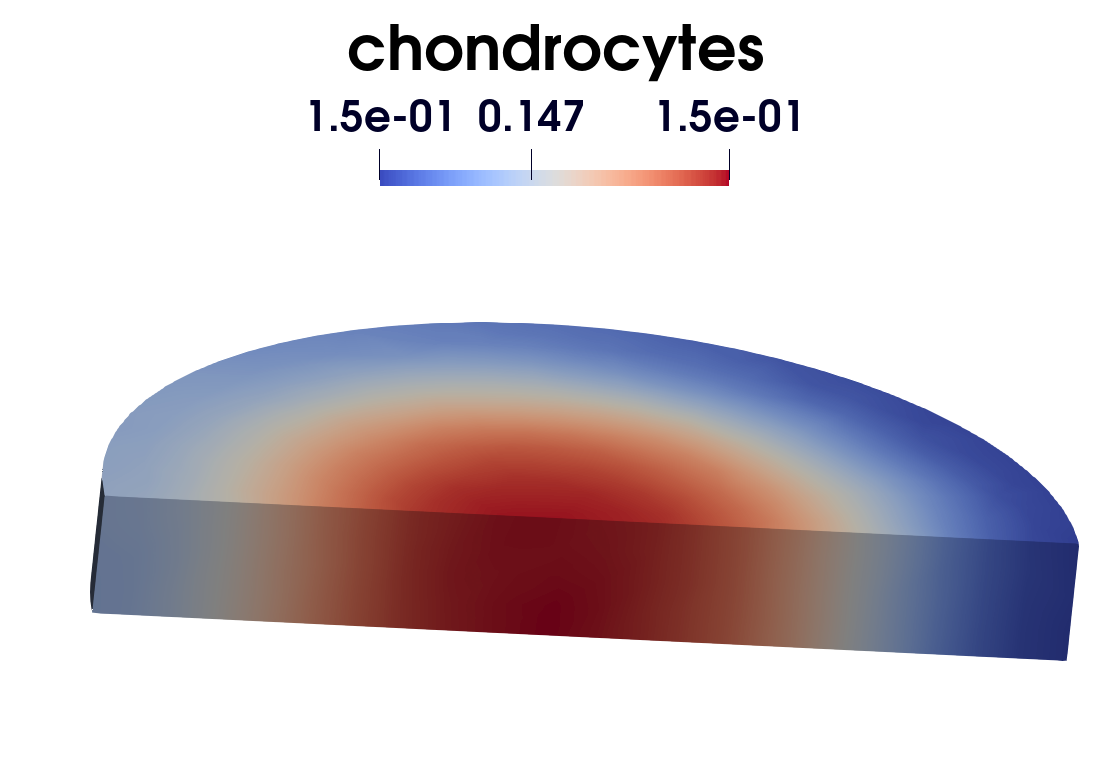}
     \end{subfigure} \quad \quad
          \begin{subfigure}[b]{0.2\textwidth}
         \centering
         \includegraphics[width=1.\textwidth]{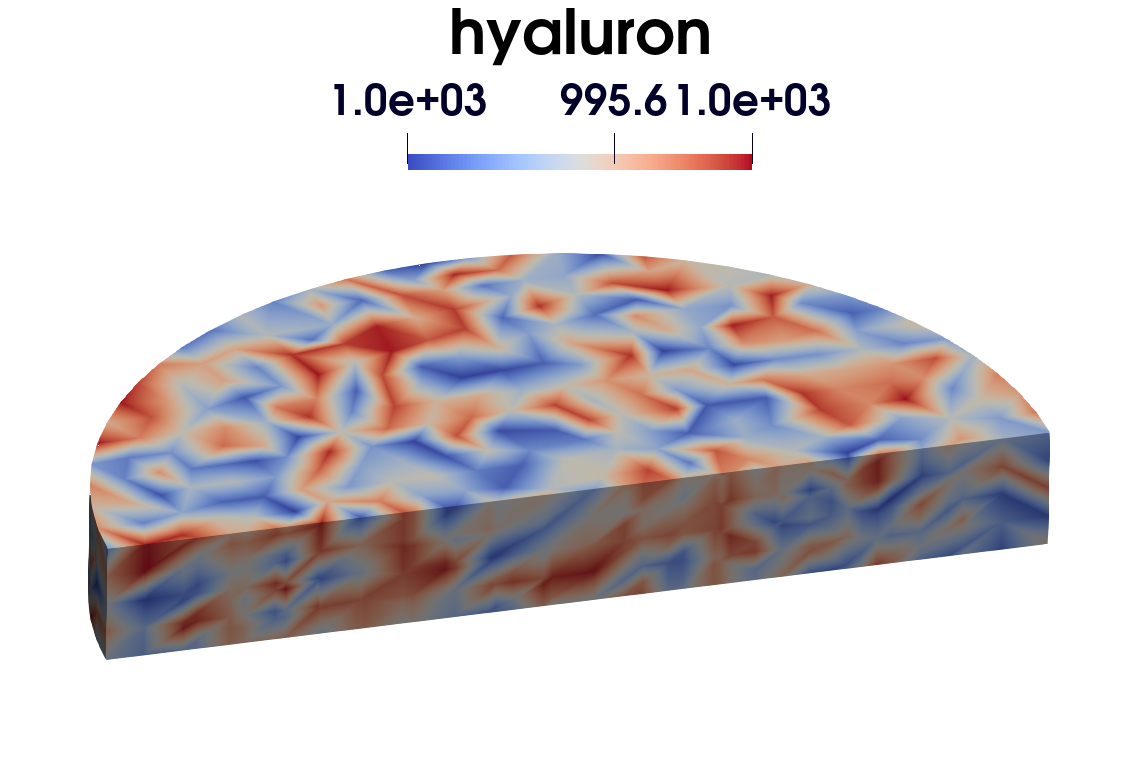}
     \end{subfigure} \quad \quad
          \begin{subfigure}[b]{0.2\textwidth}
         \centering
         \includegraphics[width=1.\textwidth]{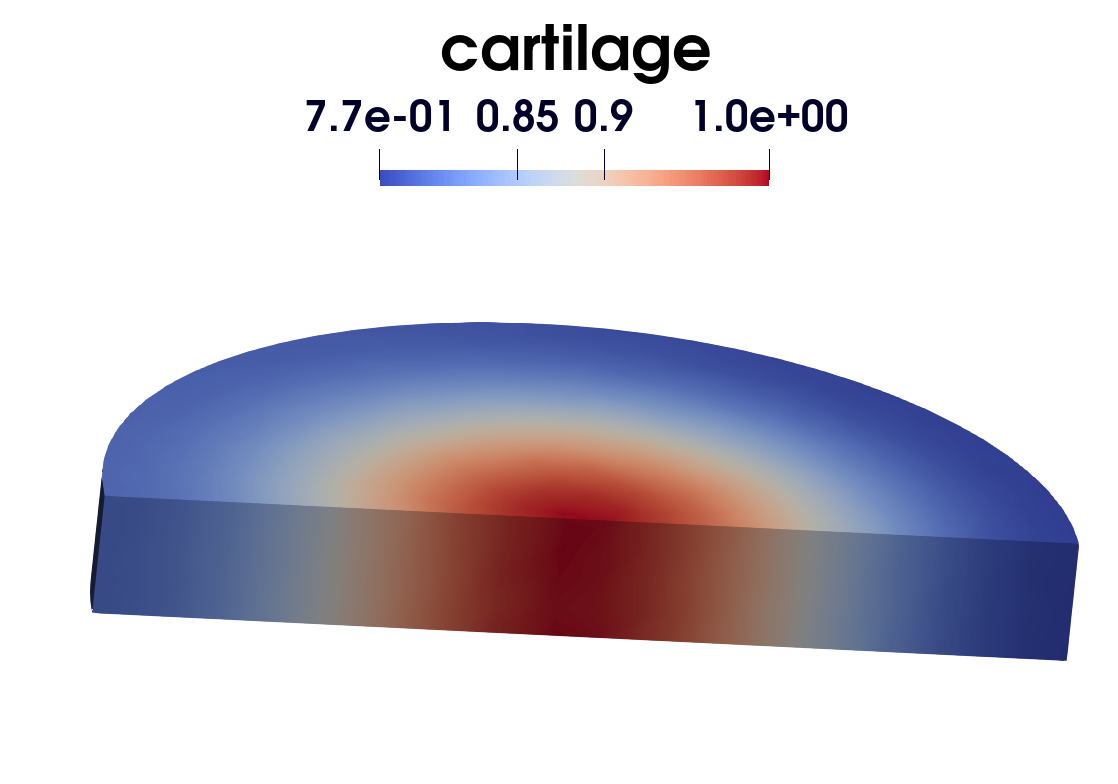}
     \end{subfigure} \\
     \vspace{-0.25cm}
         \begin{subfigure}[b]{0.2\textwidth}
         \centering
         \includegraphics[width=1.\textwidth]{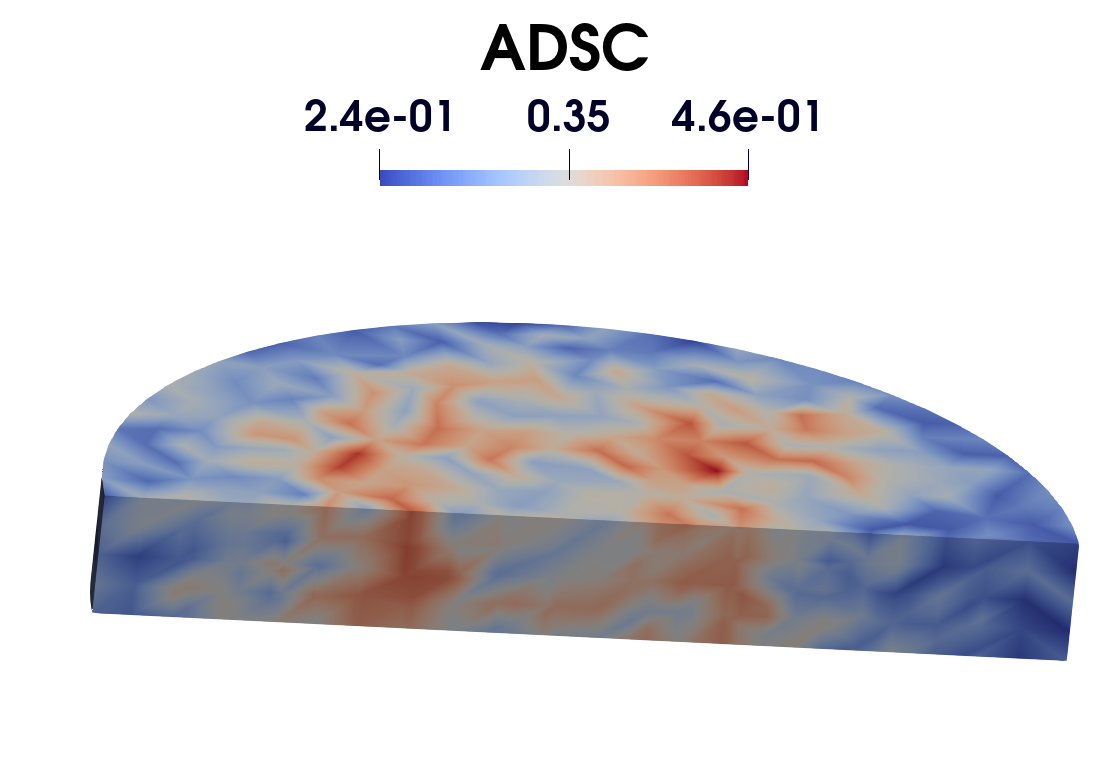} 
     \end{subfigure} \quad \quad
          \begin{subfigure}[b]{0.2\textwidth}
         \centering
         \includegraphics[width=1.\textwidth]{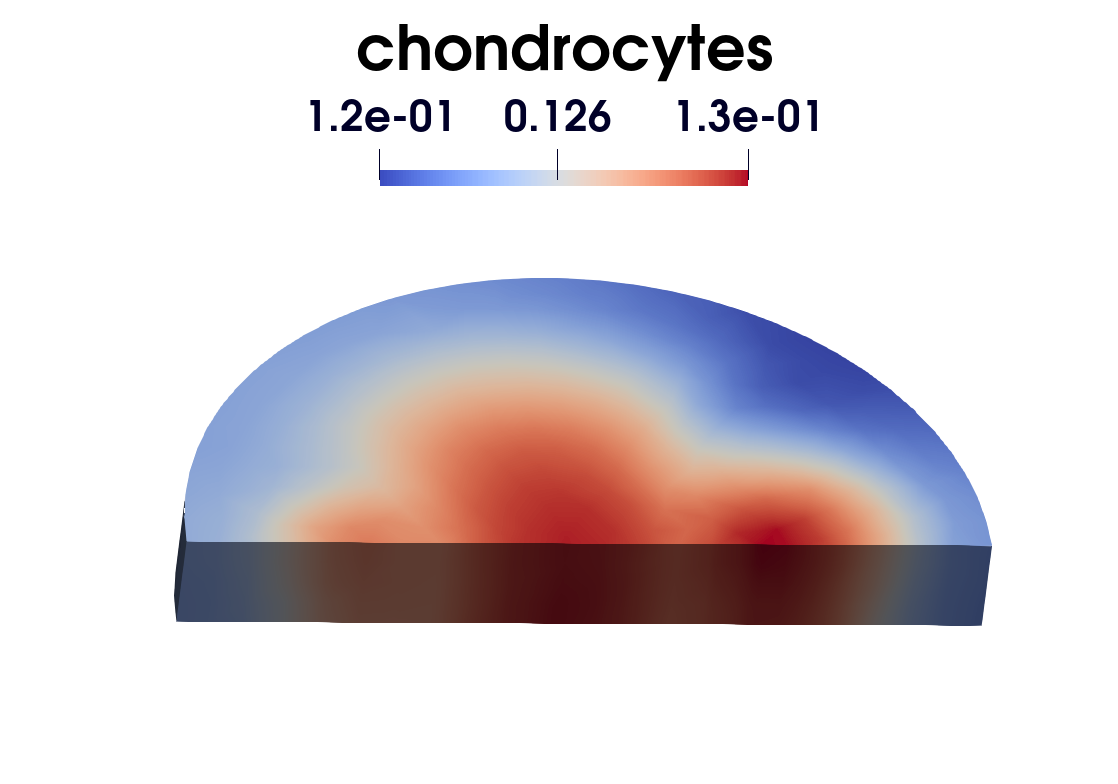}
     \end{subfigure} \quad \quad
          \begin{subfigure}[b]{0.2\textwidth}
         \centering
         \includegraphics[width=1.\textwidth]{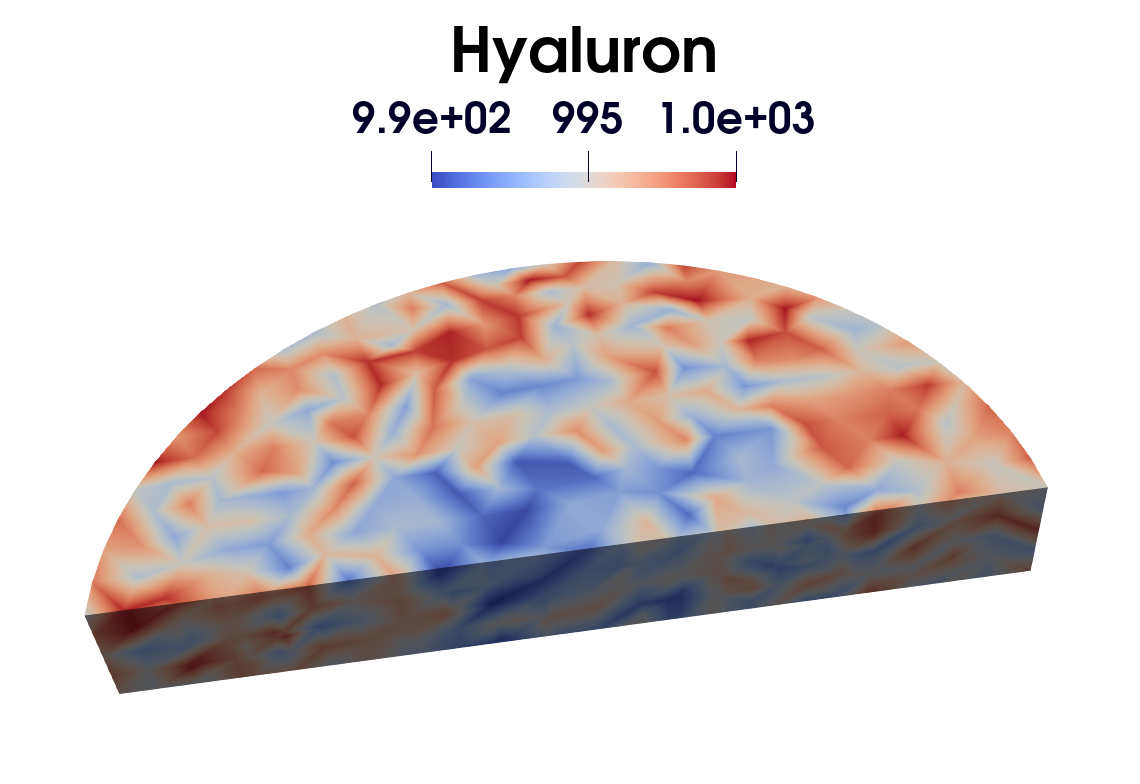}
     \end{subfigure} \quad \quad
          \begin{subfigure}[b]{0.21\textwidth}
         \centering
         \includegraphics[width=1.\textwidth]{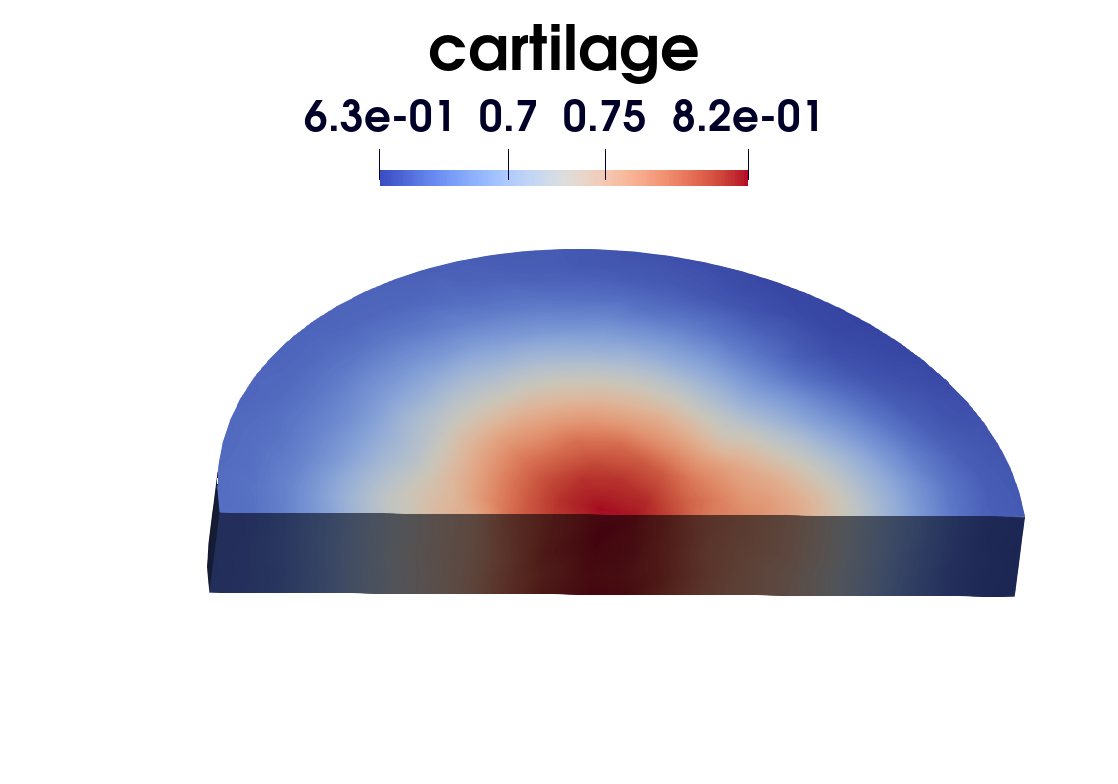}
     \end{subfigure} 
     \vspace{-0.5cm}
         \caption{After 300 iterations, with $\alpha_1=\alpha_2$ constants (up) and with $\alpha_1(\boldsymbol{\sigma}_p)=\alpha_2(\boldsymbol{\sigma}_p)$, affine functions stress-dependent (down).}
     \label{Cellresults2}
     \end{figure}

In figures \ref{Cellresults} and \ref{Cellresults2}, we plot initialization densities and stress magnitude, and  results after 300 time steps with the values of the parameters given in figure \ref{StressMap} and in table \ref{tableparam2}.\\
  \begin{minipage}{\textwidth}
 \begin{minipage}[b]{0.49\textwidth}
    \centering
      \includegraphics[width=5cm]{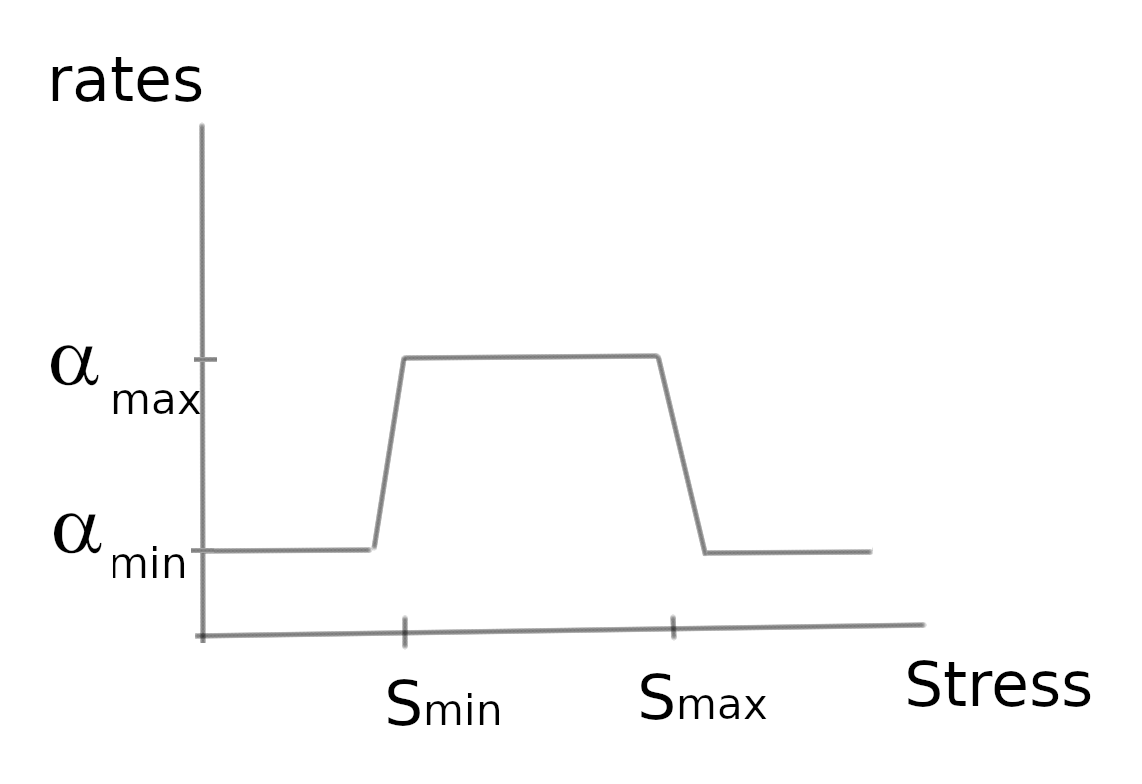}
    \captionof{figure}{Stress/rates mapping}
      \label{StressMap}
  \end{minipage}
  \hfill
     \vspace{0.25cm}
  \begin{minipage}[b]{0.49\textwidth}
    \centering
  \begin{tabular}[h]{|c|c|c|c|}\hline
        {$a_1$} &   $0.015$  &
        $ \beta$& $ 0.5$ \\
        $b_1$& $0.005$ &
        $b_2$ & $0.001$\\
      $\alpha_{min}$ & $0.05$  &
 $\alpha_{max}$ & $0.1$ \\
         $\delta_1$& $0.01 $ &
         $\gamma_1$ &$0.01$  \\
           $S_{min}$ &$1$  &
              $S_{max}$ &$3$  \\
         \hline
    \end{tabular}
    \vspace{0.6cm}
      \captionof{table}{Values of model parameters}
        \label{tableparam2}
    \end{minipage}
  \end{minipage}
  
 
\begin{minipage}{0.4\textwidth}
\vspace{-0.3cm}
\includegraphics[width=1.\textwidth]{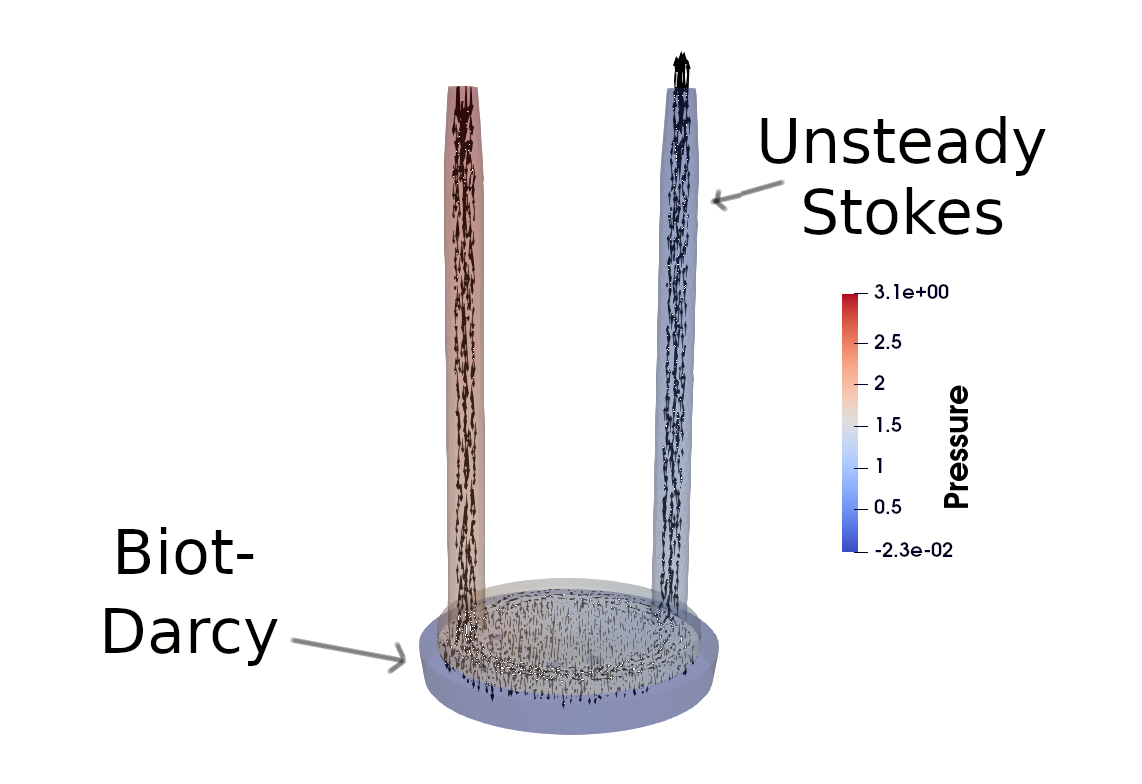}
\captionof{figure}{Pressure (MPa)  and velocity glyphs (of $\uu_f$ and $\uu_p/\Phi$) after 9 iterations. For visualization, the velocity field
on $\Omega_p$ has been magnified by a factor $5 \times 10^{7}$.}
    \label{firstplots}
\end{minipage} \;
\begin{minipage}{0.5\textwidth}
\vspace{-0.5cm}
\includegraphics[width=0.328\textwidth]{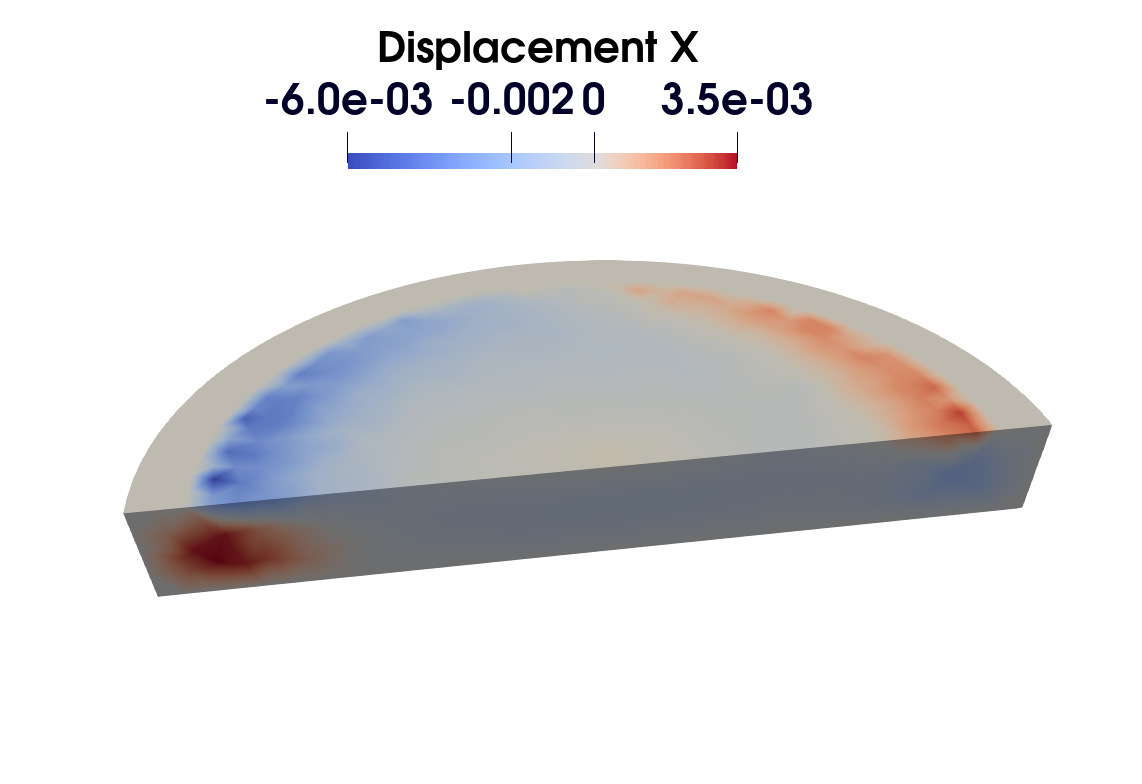}
\includegraphics[width=0.328\textwidth]{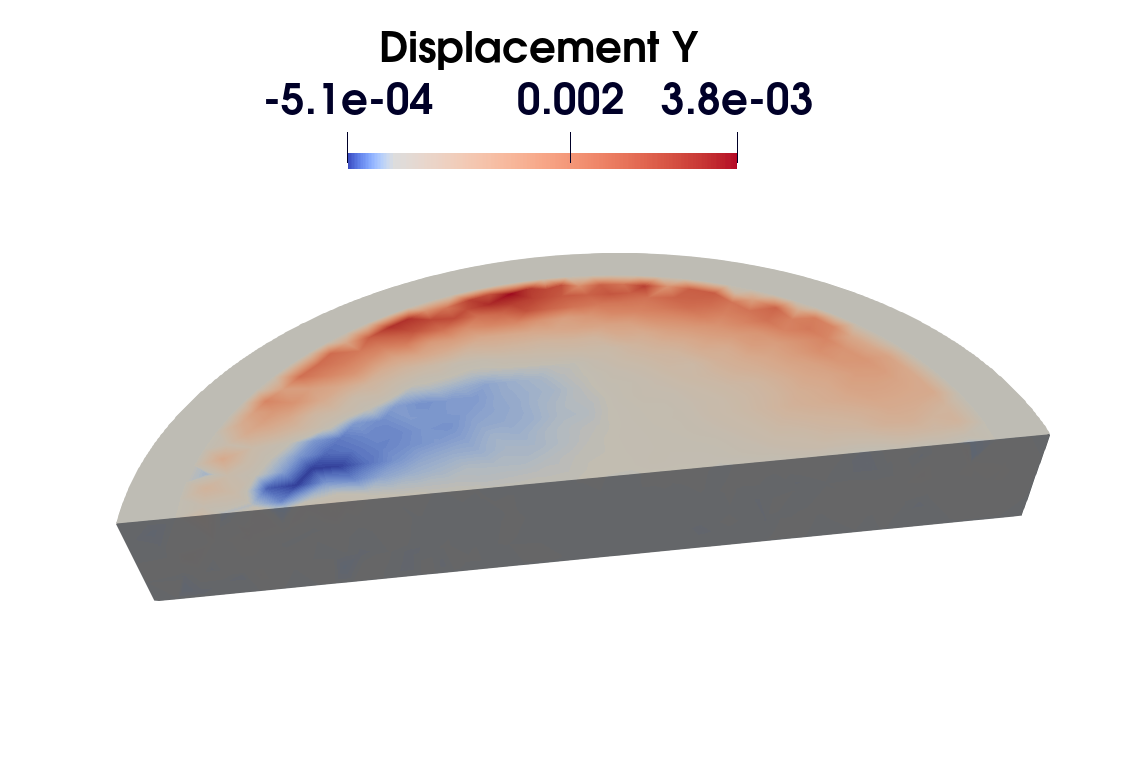}
 \includegraphics[width=0.328\textwidth]{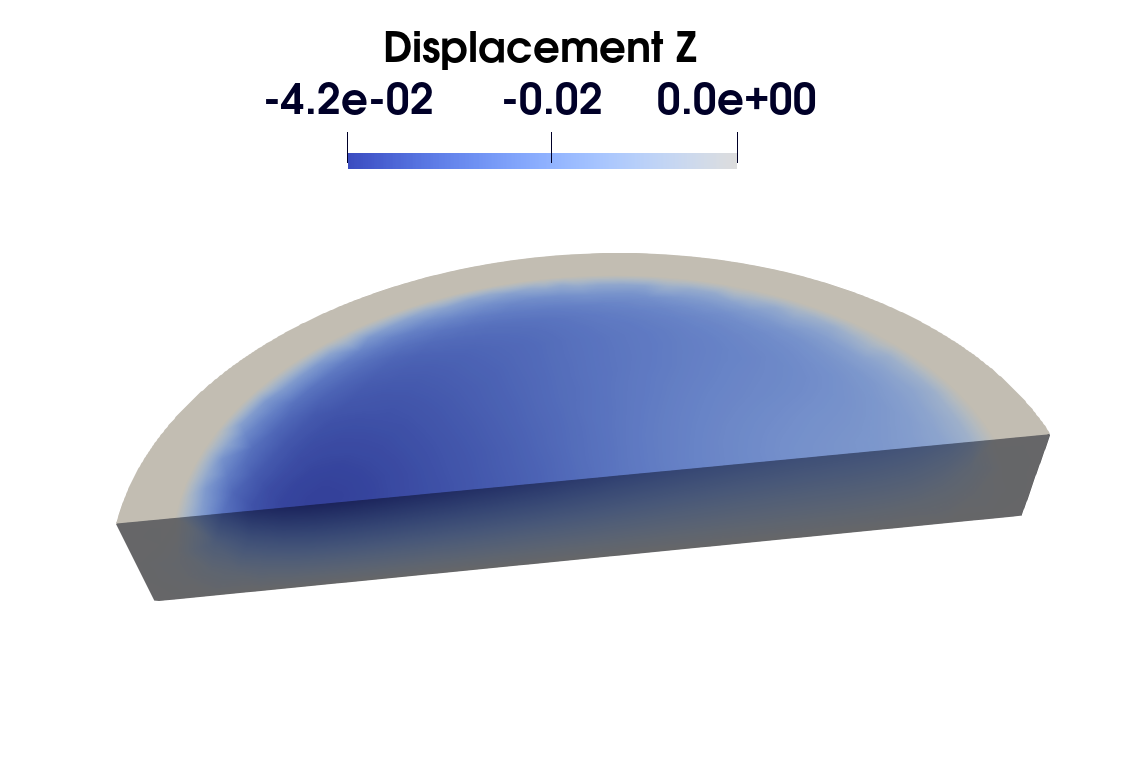}\\
\includegraphics[width=0.328\textwidth]{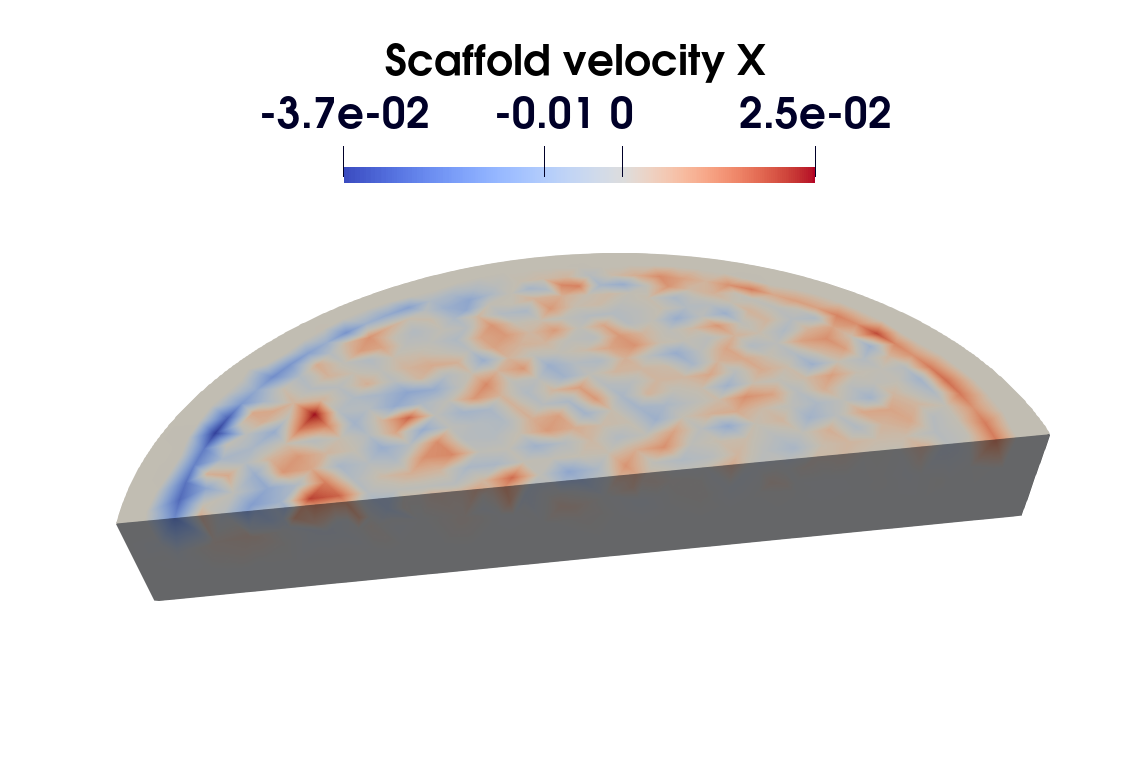}
\includegraphics[width=0.328\textwidth]{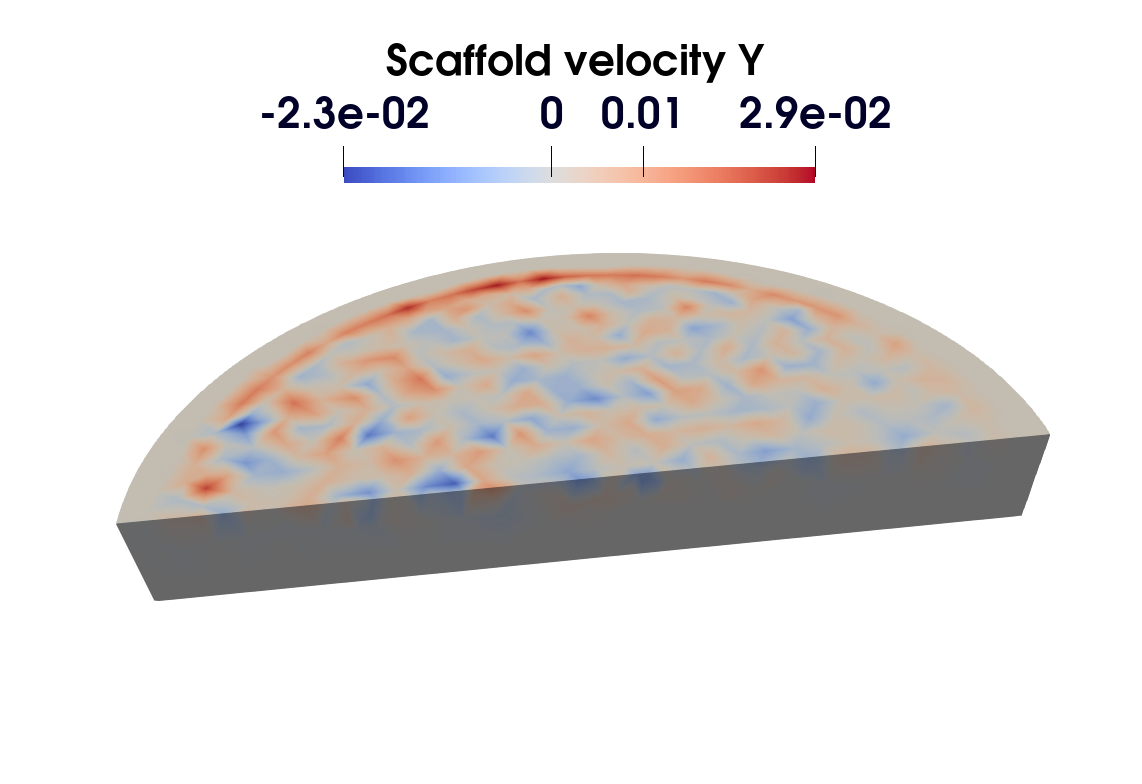}
 \includegraphics[width=0.328\textwidth]{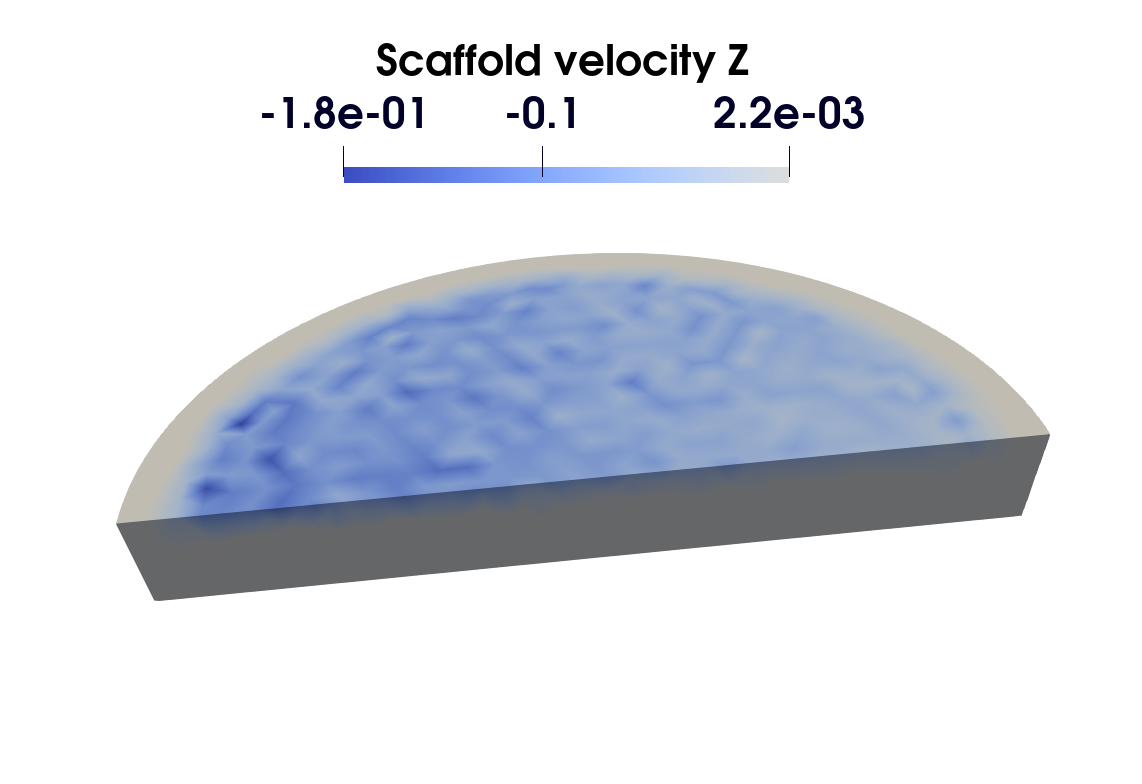}
 \vspace{-1cm}
\captionof{figure}{Displacement (mm) and velocity $\uu_p/\Phi$ (mm.s$^{-1}$) after 9 iterations (upper and lower panels, resp.)}
    \label{secondplots}
\end{minipage}

\section{Conclusions and perspectives}
Our first results with this strategy exhibit the impact of mechanical stimulus. The tissue displacements induced by the bioreactor fluid seem not to be negligible (Figure \ref{secondplots}). Moreover, Figure \ref{Cellresults2} shows that the cell densities (especially that of chondrocytes) respond to the stress magnitude. 
Finally, Figure \ref{Cellresults2} illicits that although it infers little changes over time, the chemoattractant (hyaluron) absorbed into the tissue at the start of the experiment is crucial for the ADSCs. Given how time-consuming the coupled simulations are, using a reduced order model technique \cite{NIRBSens} would be a viable option to enable calibration from in-vitro experiments. It will also be conducive when handling a more complex model explicitly accounting for the structure of the scaffold. 

\subsection{Acknowledgements}  

We are grateful to Pierre Jolivet for his assistance with PETSc and acknowledge Milena Röhrs for her study of the perfusion geometry. This work was funded by DFG within SPP 2311.

\vspace{-0.25cm}
\bibliographystyle{acm}


\end{document}